\documentclass{amsart}
\usepackage{amsfonts}

\newtheorem{theorem}{Theorem}
\theoremstyle{plain}

\newtheorem{corollary}{Corollary}

\newtheorem{proposition}{Proposition}
\newtheorem{remark}{Remark}

\numberwithin{equation}{section}
\input{tcilatex}

\begin{document}
\title[Reverses of the Generalised Triangle Inequality]{Some Reverses of the
Generalised Triangle Inequality in Complex Inner Product Spaces}
\author{Sever S. Dragomir}
\address{School of Computer Science and Mathematics\\
Victoria University of Technology\\
PO Box 14428, MCMC 8001\\
Victoria, Australia.}
\email{sever@csm.vu.edu.au}
\urladdr{http://rgmia.vu.edu.au/SSDragomirWeb.html}
\date{April 16, 2004.}

\begin{abstract}
Some reverses  for the generalised triangle inequality in complex inner
product spaces that improve the classical Diaz-Metcalf results and
applications are given.
\end{abstract}

\maketitle

\section{Introduction}

The following reverse of the generalised triangle inequality%
\begin{equation}
\cos \theta \sum_{k=1}^{n}\left\vert z_{k}\right\vert \leq \left\vert
\sum_{k=1}^{n}z_{k}\right\vert   \label{1.1}
\end{equation}%
provided the complex numbers $z_{k},$ $k\in \left\{ 1,\dots ,n\right\} $
satisfy the assumption%
\begin{equation}
a-\theta \leq \arg \left( z_{k}\right) \leq a+\theta ,\ \ \text{for any \ }%
k\in \left\{ 1,\dots ,n\right\} ,  \label{1.2}
\end{equation}%
where $a\in \mathbb{R}$ and $\theta \in \left( 0,\frac{\pi }{2}\right) $ was
first discovered by M. Petrovich in 1917, \cite{P} (see \cite[p. 492]{MPF})
and subsequently was rediscovered by other authors, including J. Karamata 
\cite[p. 300 -- 301]{K}, H.S. Wilf \cite{W}, and in an equivalent form by M.
Marden \cite{M}.

The first to consider the problem of obtaining reverses for the triangle
inequality in the more general case of Hilbert and Banach spaces were J.B.
Diaz and F.T. Metcalf \cite{DM} who showed that in an inner product space $H$
over the real or complex number field, the following reverse of the triangle
inequality holds%
\begin{equation}
r\sum_{k=1}^{n}\left\Vert x_{k}\right\Vert \leq \left\Vert
\sum_{k=1}^{n}x_{k}\right\Vert ,  \label{1.3}
\end{equation}%
provided%
\begin{equation*}
0\leq r\leq \left\Vert x_{k}\right\Vert \leq \func{Re}\left\langle
x_{k},a\right\rangle \ \ \text{for \ }k\in \left\{ 1,\dots ,n\right\} ,
\end{equation*}%
where $a\in H$ is a unit vector, i.e. $\left\Vert a\right\Vert =1.$

The case of equality holds in (\ref{1.3}) if and only if%
\begin{equation}
\sum_{k=1}^{n}x_{k}=r\left( \sum_{k=1}^{n}\left\Vert x_{k}\right\Vert
\right) a.  \label{1.4}
\end{equation}

The main purpose of this paper is to investigate the same problem of
reversing the generalised triangle inequality in complex inner product
spaces under additional assumptions for the imaginary part $\func{Im}%
\left\langle x_{k},a\right\rangle .$ A refinement of the Diaz-Metcalf result
is obtained. Applications for complex numbers are pointed out.

\section{The Case of a Unit Vector}

The following result holds.

\begin{theorem}
\label{t2.1}Let $\left( H;\left\langle \cdot ,\cdot \right\rangle \right) $
be a complex inner product space. Suppose that the vectors $x_{k}\in H,$ $%
k\in \left\{ 1,\dots ,n\right\} $ satisfy the condition%
\begin{equation}
0\leq r_{1}\left\Vert x_{k}\right\Vert \leq \func{Re}\left\langle
x_{k},e\right\rangle ,\quad 0\leq r_{2}\left\Vert x_{k}\right\Vert \leq 
\func{Im}\left\langle x_{k},e\right\rangle  \label{2.1}
\end{equation}%
for each $k\in \left\{ 1,\dots ,n\right\} ,$ where $e\in H$ is such that $%
\left\Vert e\right\Vert =1$ and $r_{1},r_{2}\geq 0.$ Then we have the
inequality%
\begin{equation}
\sqrt{r_{1}^{2}+r_{2}^{2}}\sum_{k=1}^{n}\left\Vert x_{k}\right\Vert \leq
\left\Vert \sum_{k=1}^{n}x_{k}\right\Vert ,  \label{2.2}
\end{equation}%
where equality holds if and only if%
\begin{equation}
\sum_{k=1}^{n}x_{k}=\left( r_{1}+ir_{2}\right) \left(
\sum_{k=1}^{n}\left\Vert x_{k}\right\Vert \right) e.  \label{2.3}
\end{equation}
\end{theorem}

\begin{proof}
In view of the Schwarz inequality in the complex inner product space $\left(
H;\left\langle \cdot ,\cdot \right\rangle \right) ,$ we have%
\begin{align}
\left\Vert \sum_{k=1}^{n}x_{k}\right\Vert ^{2}& =\left\Vert
\sum_{k=1}^{n}x_{k}\right\Vert ^{2}\left\Vert e\right\Vert ^{2}\geq
\left\vert \left\langle \sum_{k=1}^{n}x_{k},e\right\rangle \right\vert ^{2}
\label{2.4} \\
& =\left\vert \left\langle \sum_{k=1}^{n}x_{k},e\right\rangle \right\vert
^{2}  \notag \\
& =\left\vert \sum_{k=1}^{n}\func{Re}\left\langle x_{k},e\right\rangle
+i\left( \sum_{k=1}^{n}\func{Im}\left\langle x_{k},e\right\rangle \right)
\right\vert ^{2}  \notag \\
& =\left( \sum_{k=1}^{n}\func{Re}\left\langle x_{k},e\right\rangle \right)
^{2}+\left( \sum_{k=1}^{n}\func{Im}\left\langle x_{k},e\right\rangle \right)
^{2}.  \notag
\end{align}%
Now, by hypothesis (\ref{2.1})%
\begin{equation}
\left( \sum_{k=1}^{n}\func{Re}\left\langle x_{k},e\right\rangle \right)
^{2}\geq r_{1}^{2}\left( \sum_{k=1}^{n}\left\Vert x_{k}\right\Vert \right)
^{2}  \label{2.5}
\end{equation}%
and%
\begin{equation}
\left( \sum_{k=1}^{n}\func{Im}\left\langle x_{k},e\right\rangle \right)
^{2}\geq r_{2}^{2}\left( \sum_{k=1}^{n}\left\Vert x_{k}\right\Vert \right)
^{2}.  \label{2.6}
\end{equation}%
If we add (\ref{2.5}) and (\ref{2.6}) and use (\ref{2.4}), then we deduce
the desired inequality (\ref{2.2}).

Now, if (\ref{2.3}) holds, then%
\begin{equation*}
\left\Vert \sum_{k=1}^{n}x_{k}\right\Vert =\left\vert
r_{1}+ir_{2}\right\vert \left( \sum_{k=1}^{n}\left\Vert x_{k}\right\Vert
\right) \left\Vert e\right\Vert =\sqrt{r_{1}^{2}+r_{2}^{2}}%
\sum_{k=1}^{n}\left\Vert x_{k}\right\Vert 
\end{equation*}%
and the case of equality is valid in (\ref{2.2}).

Before we prove the reverse implication, let us observe that for $x\in H$
and $e\in H,$ $\left\Vert e\right\Vert =1,$ the following identity is true%
\begin{equation*}
\left\Vert x-\left\langle x,e\right\rangle e\right\Vert ^{2}=\left\Vert
x\right\Vert ^{2}-\left\vert \left\langle x,e\right\rangle \right\vert ^{2},
\end{equation*}%
therefore $\left\Vert x\right\Vert =\left\vert \left\langle x,e\right\rangle
\right\vert $ if and only if $x=\left\langle x,e\right\rangle e.$

If we assume that equality holds in (\ref{2.2}), then the case of equality
must hold in all the inequalities required in the argument used to prove the
inequality (\ref{1.2}), and we may state that%
\begin{equation}
\left\Vert \sum_{k=1}^{n}x_{k}\right\Vert =\left\vert \left\langle
\sum_{k=1}^{n}x_{k},e\right\rangle \right\vert ,  \label{2.7}
\end{equation}%
and 
\begin{equation}
r_{1}\left\Vert x_{k}\right\Vert =\func{Re}\left\langle x_{k},e\right\rangle
,\quad r_{2}\left\Vert x_{k}\right\Vert =\func{Im}\left\langle
x_{k},e\right\rangle  \label{2.8}
\end{equation}%
for each $k\in \left\{ 1,\dots ,n\right\} .$

From (\ref{2.7}) we deduce%
\begin{equation}
\sum_{k=1}^{n}x_{k}=\left\langle \sum_{k=1}^{n}x_{k},e\right\rangle e
\label{2.9}
\end{equation}%
and from (\ref{2.8}), by multiplying the second equation with $i$ and
summing both equations over $k$ from $1$ to $n,$ we deduce%
\begin{equation}
\left( r_{1}+ir_{2}\right) \sum_{k=1}^{n}\left\Vert x_{k}\right\Vert
=\left\langle \sum_{k=1}^{n}x_{k},e\right\rangle .  \label{2.10}
\end{equation}%
Finally, by (\ref{2.10}) and (\ref{2.9}), we get the desired equality (\ref%
{2.3}).
\end{proof}

The following corollary is of interest.

\begin{corollary}
\label{c2.2}Let $e$ a unit vector in the complex inner product space $\left(
H;\left\langle \cdot ,\cdot \right\rangle \right) $ and $\rho _{1},\rho
_{2}\in \left( 0,1\right) .$ If $x_{k}\in H,$ $k\in \left\{ 1,\dots
,n\right\} $ are such that%
\begin{equation}
\left\Vert x_{k}-e\right\Vert \leq \rho _{1},\ \ \ \left\Vert
x_{k}-ie\right\Vert \leq \rho _{2}\ \ \ \text{for each \ }k\in \left\{
1,\dots ,n\right\} ,  \label{2.11}
\end{equation}%
then we have the inequality%
\begin{equation}
\sqrt{2-\rho _{1}^{2}-\rho _{2}^{2}}\sum_{k=1}^{n}\left\Vert
x_{k}\right\Vert \leq \left\Vert \sum_{k=1}^{n}x_{k}\right\Vert ,
\label{2.12}
\end{equation}%
with equality if and only if%
\begin{equation}
\sum_{k=1}^{n}x_{k}=\left( \sqrt{1-\rho _{1}^{2}}+i\sqrt{1-\rho _{2}^{2}}%
\right) \left( \sum_{k=1}^{n}\left\Vert x_{k}\right\Vert \right) e.
\label{2.13}
\end{equation}
\end{corollary}

\begin{proof}
From the first inequality in (\ref{2.11}) we deduce, by taking the square,
that%
\begin{equation*}
\left\Vert x_{k}\right\Vert ^{2}+1-\rho _{1}^{2}\leq 2\func{Re}\left\langle
x_{k},e\right\rangle ,
\end{equation*}%
implying%
\begin{equation}
\frac{\left\Vert x_{k}\right\Vert ^{2}}{\sqrt{1-\rho _{1}^{2}}}+\sqrt{1-\rho
_{1}^{2}}\leq \frac{2\func{Re}\left\langle x_{k},e\right\rangle }{\sqrt{%
1-\rho _{1}^{2}}}  \label{2.14}
\end{equation}%
for each $k\in \left\{ 1,\dots ,n\right\} .$

Since, obviously,%
\begin{equation}
2\left\Vert x_{k}\right\Vert \leq \frac{\left\Vert x_{k}\right\Vert ^{2}}{%
\sqrt{1-\rho _{1}^{2}}}+\sqrt{1-\rho _{1}^{2}},\ \text{\ }k\in \left\{
1,\dots ,n\right\} ,  \label{2.15}
\end{equation}%
hence, by (\ref{2.14}) and (\ref{2.5}),%
\begin{equation}
0\leq \sqrt{1-\rho _{1}^{2}}\left\Vert x_{k}\right\Vert \leq \func{Re}%
\left\langle x_{k},e\right\rangle   \label{2.16}
\end{equation}%
for each $k\in \left\{ 1,\dots ,n\right\} .$

From the second inequality in (\ref{2.11}) we deduce%
\begin{equation*}
0\leq \sqrt{1-\rho _{2}^{2}}\left\Vert x_{k}\right\Vert \leq \func{Re}%
\left\langle x_{k},ie\right\rangle 
\end{equation*}%
for each $k\in \left\{ 1,\dots ,n\right\} .$ Since%
\begin{equation*}
\func{Re}\left\langle x_{k},ie\right\rangle =\func{Im}\left\langle
x_{k},e\right\rangle ,
\end{equation*}%
hence%
\begin{equation}
0\leq \sqrt{1-\rho _{2}^{2}}\left\Vert x_{k}\right\Vert \leq \func{Im}%
\left\langle x_{k},e\right\rangle   \label{2.17}
\end{equation}%
for each $k\in \left\{ 1,\dots ,n\right\} .$

Now, observe from (\ref{2.16}) and (\ref{2.17}), that the condition (\ref%
{2.1}) of Theorem \ref{t2.1} is satisfied for $r_{1}=\sqrt{1-\rho _{1}^{2}},$
$r_{2}=\sqrt{1-\rho _{2}^{2}}\in \left( 0,1\right) ,$ and thus the corollary
is proved.
\end{proof}

The following corollary may be stated as well.

\begin{corollary}
\label{c2.3}Let $e$ be a unit vector in the complex inner product space $%
\left( H;\left\langle \cdot ,\cdot \right\rangle \right) $ and $M_{1}\geq
m_{1}>0,$ $M_{2}\geq m_{2}>0.$ If $x_{k}\in H,$ $k\in \left\{ 1,\dots
,n\right\} $ are such that either%
\begin{equation}
\func{Re}\left\langle M_{1}e-x_{k},x_{k}-m_{1}e\right\rangle \geq 0,\ \ \ 
\text{ \ }\func{Re}\left\langle M_{2}ie-x_{k},x_{k}-m_{2}ie\right\rangle
\geq 0  \label{2.18}
\end{equation}%
or, equivalently,%
\begin{align}
\left\Vert x_{k}-\frac{M_{1}+m_{1}}{2}e\right\Vert & \leq \frac{1}{2}\left(
M_{1}-m_{1}\right) ,  \label{2.19} \\
\left\Vert x_{k}-\frac{M_{2}+m_{2}}{2}ie\right\Vert & \leq \frac{1}{2}\left(
M_{2}-m_{2}\right) ,  \notag
\end{align}%
for each $k\in \left\{ 1,\dots ,n\right\} ,$ then we have the inequality%
\begin{equation}
2\left[ \frac{m_{1}M_{1}}{\left( M_{1}+m_{1}\right) ^{2}}+\frac{m_{2}M_{2}}{%
\left( M_{2}+m_{2}\right) ^{2}}\right] ^{1/2}\sum_{k=1}^{n}\left\Vert
x_{k}\right\Vert \leq \left\Vert \sum_{k=1}^{n}x_{k}\right\Vert .
\label{2.20}
\end{equation}%
The equality holds in (\ref{2.20}) if and only if%
\begin{equation}
\sum_{k=1}^{n}x_{k}=2\left( \frac{\sqrt{m_{1}M_{1}}}{M_{1}+m_{1}}+i\frac{%
\sqrt{m_{2}M_{2}}}{M_{2}+m_{2}}\right) \left( \sum_{k=1}^{n}\left\Vert
x_{k}\right\Vert \right) e.  \label{2.21}
\end{equation}
\end{corollary}

\begin{proof}
Firstly, remark that, for $x,z,Z\in H,$ the following statements are
equivalent.

\begin{enumerate}
\item[(i)] $\func{Re}\left\langle Z-x,x-z\right\rangle \geq 0$

and

\item[(ii)] $\left\Vert x-\frac{Z+z}{2}\right\Vert \leq \frac{1}{2}%
\left\Vert Z-z\right\Vert .$
\end{enumerate}

Using this fact, we may simply realise that (\ref{2.18}) and (\ref{2.19})
are equivalent.

Now, from the first inequality in (\ref{2.18}), we get%
\begin{equation*}
\left\Vert x_{k}\right\Vert ^{2}+m_{1}M_{1}\leq \left( M_{1}+m_{1}\right) 
\func{Re}\left\langle x_{k},e\right\rangle
\end{equation*}%
implying%
\begin{equation}
\frac{\left\Vert x_{k}\right\Vert ^{2}}{\sqrt{m_{1}M_{1}}}+\sqrt{m_{1}M_{1}}%
\leq \frac{M_{1}+m_{1}}{\sqrt{m_{1}M_{1}}}\func{Re}\left\langle
x_{k},e\right\rangle  \label{2.22}
\end{equation}%
for each $k\in \left\{ 1,\dots ,n\right\} .$

Since, obviously,%
\begin{equation}
2\left\Vert x_{k}\right\Vert \leq \frac{\left\Vert x_{k}\right\Vert ^{2}}{%
\sqrt{m_{1}M_{1}}}+\sqrt{m_{1}M_{1}},  \label{2.23}
\end{equation}%
hence, by (\ref{2.22}) and (\ref{2.23})%
\begin{equation}
0\leq \frac{2\sqrt{m_{1}M_{1}}}{M_{1}+m_{1}}\left\Vert x_{k}\right\Vert \leq 
\func{Re}\left\langle x_{k},e\right\rangle  \label{2.24}
\end{equation}%
for each $k\in \left\{ 1,\dots ,n\right\} .$

Now, the proof follows the same path as the one of Corollary \ref{c2.2} and
we omit the details.
\end{proof}

\section{The Case of $m$ Orthornormal Vectors}

In \cite{DM}, the authors have proved the following reverse of the
generalised triangle inequality in terms of orthornormal vectors.

\begin{theorem}
\label{t3.1}Let $e_{1},\dots ,e_{m}$ be orthornormal vectors in $\left(
H;\left\langle \cdot ,\cdot \right\rangle \right) $, i.e., we recall that $%
\left\langle e_{i},e_{j}\right\rangle =0$ if $i\neq j$ and $\left\Vert
e_{i}\right\Vert =1,$ $i,j\in \left\{ 1,\dots ,m\right\} .$ Suppose that the
vectors $x_{1},\dots ,x_{n}\in H$ satisfy%
\begin{equation}
0\leq r_{k}\left\Vert x_{j}\right\Vert \leq \func{Re}\left\langle
x_{j},e_{k}\right\rangle ,\qquad j\in \left\{ 1,\dots ,n\right\} ,\ k\in
\left\{ 1,\dots ,m\right\} .  \label{3.1}
\end{equation}%
Then%
\begin{equation}
\left( \sum_{k=1}^{m}r_{k}^{2}\right) ^{\frac{1}{2}}\sum_{j=1}^{n}\left\Vert
x_{j}\right\Vert \leq \left\Vert \sum_{j=1}^{n}x_{j}\right\Vert ,
\label{3.2}
\end{equation}%
where equality holds if and only if%
\begin{equation}
\sum_{j=1}^{n}x_{j}=\left( \sum_{j=1}^{n}\left\Vert x_{j}\right\Vert \right)
\sum_{k=1}^{m}r_{k}e_{k}.  \label{3.3}
\end{equation}
\end{theorem}

If the space $\left( H;\left\langle \cdot ,\cdot \right\rangle \right) $ is
complex and more information is available for the imaginary part, then the
following result may be stated as well.

\begin{theorem}
\label{t3.2}Let $e_{1},\dots ,e_{m}\in H$ be an orthornormal family of
vectors in the complex inner product space $H.$ If the vectors $x_{1},\dots
,x_{n}\in H$ satisfy the conditions%
\begin{equation}
0\leq r_{k}\left\Vert x_{j}\right\Vert \leq \func{Re}\left\langle
x_{j},e_{k}\right\rangle ,\qquad 0\leq \rho _{k}\left\Vert x_{j}\right\Vert
\leq \func{Im}\left\langle x_{j},e_{k}\right\rangle  \label{3.4}
\end{equation}%
for each $j\in \left\{ 1,\dots ,n\right\} $ and $k\in \left\{ 1,\dots
,m\right\} ,$ then we have the following reverse of the generalised triangle
inequality;%
\begin{equation}
\left[ \sum_{k=1}^{m}\left( r_{k}^{2}+\rho _{k}^{2}\right) \right] ^{\frac{1%
}{2}}\sum_{j=1}^{n}\left\Vert x_{j}\right\Vert \leq \left\Vert
\sum_{j=1}^{n}x_{j}\right\Vert .  \label{3.5}
\end{equation}%
The equality holds in (\ref{3.5}) if and only if%
\begin{equation}
\sum_{j=1}^{n}x_{j}=\left( \sum_{j=1}^{n}\left\Vert x_{j}\right\Vert \right)
\sum_{k=1}^{m}\left( r_{k}+i\rho _{k}\right) e_{k}.  \label{3.6}
\end{equation}
\end{theorem}

\begin{proof}
Before we prove the theorem, let us recall that, if $x\in H$ and $%
e_{1},\dots ,e_{m}$ are orthogonal vectors, then the following identity
holds true:%
\begin{equation}
\left\Vert x-\sum_{k=1}^{m}\left\langle x,e_{k}\right\rangle
e_{k}\right\Vert ^{2}=\left\Vert x\right\Vert ^{2}-\sum_{k=1}^{n}\left\vert
\left\langle x,e_{k}\right\rangle \right\vert ^{2}.  \label{3.7}
\end{equation}%
As a consequence of this identity, we note the \textit{Bessel inequality}%
\begin{equation}
\sum_{k=1}^{m}\left\vert \left\langle x,e_{k}\right\rangle \right\vert
^{2}\leq \left\Vert x\right\Vert ^{2},x\in H.  \label{3.8}
\end{equation}%
The case of equality holds in (\ref{3.8}) if and only if (see (\ref{3.7}))%
\begin{equation}
x=\sum_{k=1}^{m}\left\langle x,e_{k}\right\rangle e_{k}.  \label{3.9}
\end{equation}%
Applying Bessel's inequality for $x=\sum_{j=1}^{n}x_{j},$ we have%
\begin{align}
\left\Vert \sum_{j=1}^{n}x_{j}\right\Vert ^{2}& \geq
\sum_{k=1}^{m}\left\vert \left\langle \sum_{j=1}^{n}x_{j},e_{k}\right\rangle
\right\vert ^{2}=\sum_{k=1}^{m}\left\vert \sum_{j=1}^{n}\left\langle
x_{j},e_{k}\right\rangle \right\vert ^{2}  \label{3.10} \\
& =\sum_{k=1}^{m}\left\vert \left( \sum_{j=1}^{n}\func{Re}\left\langle
x_{j},e_{k}\right\rangle \right) +i\left( \sum_{j=1}^{n}\func{Im}%
\left\langle x_{j},e_{k}\right\rangle \right) \right\vert ^{2}  \notag \\
& =\sum_{k=1}^{m}\left[ \left( \sum_{j=1}^{n}\func{Re}\left\langle
x_{j},e_{k}\right\rangle \right) ^{2}+\left( \sum_{j=1}^{n}\func{Im}%
\left\langle x_{j},e_{k}\right\rangle \right) ^{2}\right] .  \notag
\end{align}%
Now, by the hypothesis (\ref{3.4}) we have%
\begin{equation}
\left( \sum_{j=1}^{n}\func{Re}\left\langle x_{j},e_{k}\right\rangle \right)
^{2}\geq r_{k}^{2}\left( \sum_{j=1}^{n}\left\Vert x_{j}\right\Vert \right)
^{2}  \label{3.11}
\end{equation}%
and%
\begin{equation}
\left( \sum_{j=1}^{n}\func{Im}\left\langle x_{j},e_{k}\right\rangle \right)
^{2}\geq \rho _{k}^{2}\left( \sum_{j=1}^{n}\left\Vert x_{j}\right\Vert
\right) ^{2}.  \label{3.12}
\end{equation}%
Further, on making use of (\ref{3.10}) -- (\ref{3.12}), we deduce%
\begin{align*}
\left\Vert \sum_{j=1}^{n}x_{j}\right\Vert ^{2}& \geq \sum_{k=1}^{m}\left[
r_{k}^{2}\left( \sum_{j=1}^{n}\left\Vert x_{j}\right\Vert \right) ^{2}+\rho
_{k}^{2}\left( \sum_{j=1}^{n}\left\Vert x_{j}\right\Vert \right) ^{2}\right] 
\\
& =\left( \sum_{j=1}^{n}\left\Vert x_{j}\right\Vert \right)
^{2}\sum_{k=1}^{m}\left( r_{k}^{2}+\rho _{k}^{2}\right) ,
\end{align*}%
which is clearly equivalent to (\ref{3.5}).

Now, if (\ref{3.6}) holds, then%
\begin{align*}
\left\Vert \sum_{j=1}^{n}x_{j}\right\Vert ^{2}& =\left(
\sum_{j=1}^{n}\left\Vert x_{j}\right\Vert \right) ^{2}\left\Vert
\sum_{k=1}^{m}\left( r_{k}+i\rho _{k}\right) e_{k}\right\Vert ^{2} \\
& =\left( \sum_{j=1}^{n}\left\Vert x_{j}\right\Vert \right)
^{2}\sum_{k=1}^{m}\left\vert r_{k}+i\rho _{k}\right\vert ^{2} \\
& =\left( \sum_{j=1}^{n}\left\Vert x_{j}\right\Vert \right)
^{2}\sum_{k=1}^{m}\left( r_{k}^{2}+\rho _{k}^{2}\right) ,
\end{align*}%
and the case of equality holds in (\ref{3.5}).

Conversely, if the equality holds in (\ref{3.5}), then it must hold in all
the inequalities used to prove (\ref{3.5}) and therefore we must have%
\begin{equation}
\left\Vert \sum_{j=1}^{n}x_{j}\right\Vert ^{2}=\sum_{k=1}^{m}\left\vert
\sum_{j=1}^{n}\left\langle x_{j},e_{k}\right\rangle \right\vert ^{2}
\label{3.13}
\end{equation}%
and%
\begin{equation}
r_{k}\left\Vert x_{j}\right\Vert =\func{Re}\left\langle
x_{j},e_{k}\right\rangle ,\qquad \rho _{k}\left\Vert x_{j}\right\Vert =\func{%
Im}\left\langle x_{j},e_{k}\right\rangle   \label{3.14}
\end{equation}%
for each $j\in \left\{ 1,\dots ,n\right\} $ and $k\in \left\{ 1,\dots
,m\right\} .$

Using the identity (\ref{3.7}), we deduce from (\ref{3.13}) that%
\begin{equation}
\sum_{j=1}^{n}x_{j}=\sum_{k=1}^{m}\left\langle
\sum_{j=1}^{n}x_{j},e_{k}\right\rangle e_{k}.  \label{3.15}
\end{equation}%
Multiplying the second equality in (\ref{3.14}) with the imaginary unit $i$
and summing the equality over $j$ from $1$ to $n,$ we deduce%
\begin{equation}
\left( r_{k}+i\rho _{k}\right) \sum_{j=1}^{n}\left\Vert x_{j}\right\Vert
=\left\langle \sum_{j=1}^{n}x_{j},e_{k}\right\rangle  \label{3.16}
\end{equation}%
for each $k\in \left\{ 1,\dots ,n\right\} .$

Finally, utilising (\ref{3.15}) and (\ref{3.16}), we deduce (\ref{3.6}) and
the theorem is proved.
\end{proof}

The following corollaries are of interest.

\begin{corollary}
\label{c3.2}Let $e_{1},\dots ,e_{m}$ be orthornormal vectors in the complex
inner product space $\left( H;\left\langle \cdot ,\cdot \right\rangle
\right) $ and $\rho _{k},\eta _{k}\in \left( 0,1\right) ,$ $k\in \left\{
1,\dots ,n\right\} .$ If $x_{1},\dots ,x_{n}\in H$ are such that%
\begin{equation*}
\left\Vert x_{j}-e_{k}\right\Vert \leq \rho _{k},\qquad \left\Vert
x_{j}-ie_{k}\right\Vert \leq \eta _{k}
\end{equation*}%
for each $j\in \left\{ 1,\dots ,n\right\} $ and $k\in \left\{ 1,\dots
,m\right\} ,$ then we have the inequality%
\begin{equation}
\left[ \sum_{k=1}^{m}\left( 2-\rho _{k}^{2}-\eta _{k}^{2}\right) \right] ^{%
\frac{1}{2}}\sum_{j=1}^{n}\left\Vert x_{j}\right\Vert \leq \left\Vert
\sum_{j=1}^{n}x_{j}\right\Vert .  \label{3.17}
\end{equation}%
The case of equality holds in (\ref{3.17}) if and only if%
\begin{equation}
\sum_{j=1}^{n}x_{j}=\left( \sum_{j=1}^{n}\left\Vert x_{j}\right\Vert \right)
\sum_{k=1}^{m}\left( \sqrt{1-\rho _{k}^{2}}+i\sqrt{1-\eta _{k}^{2}}\right)
e_{k}.  \label{3.18}
\end{equation}
\end{corollary}

The proof employs Theorem \ref{t3.2} and is similar to the one from
Corollary \ref{c2.2}. We omit the details.

\begin{corollary}
\label{c3.3}Let $e_{1},\dots ,e_{m}$ be as in Corollary \ref{c3.2} and $%
M_{k}\geq m_{k}>0,$ $N_{k}\geq n_{k}>0,$ $k\in \left\{ 1,\dots ,m\right\} .$
If $x_{1},\dots ,x_{n}\in H$ are such that either%
\begin{equation*}
\func{Re}\left\langle M_{k}e_{k}-x_{j},x_{j}-m_{k}e_{k}\right\rangle \geq
0,\ \ \func{Re}\left\langle N_{k}ie_{k}-x_{j},x_{j}-n_{k}ie_{k}\right\rangle
\geq 0
\end{equation*}%
or, equivalently,%
\begin{align*}
\left\Vert x_{j}-\frac{M_{k}+m_{k}}{2}e_{k}\right\Vert & \leq \frac{1}{2}%
\left( M_{k}-m_{k}\right) ,\  \\
\left\Vert x_{j}-\frac{N_{k}+n_{k}}{2}ie_{k}\right\Vert & \leq \frac{1}{2}%
\left( N_{k}-n_{k}\right) 
\end{align*}%
for each $j\in \left\{ 1,\dots ,n\right\} $ and $k\in \left\{ 1,\dots
,m\right\} ,$ then we have the inequality%
\begin{equation}
2\left\{ \sum_{k=1}^{m}\left[ \frac{m_{k}M_{k}}{\left( M_{k}+m_{k}\right)
^{2}}+\frac{n_{k}N_{k}}{\left( N_{k}+n_{k}\right) ^{2}}\right] \right\} ^{%
\frac{1}{2}}\sum_{j=1}^{n}\left\Vert x_{j}\right\Vert \leq \left\Vert
\sum_{j=1}^{n}x_{j}\right\Vert .  \label{3.19}
\end{equation}%
The case of equality holds in (\ref{3.19}) if and only if%
\begin{equation}
\sum_{j=1}^{n}x_{j}=2\left( \sum_{j=1}^{n}\left\Vert x_{j}\right\Vert
\right) \sum_{k=1}^{m}\left( \frac{\sqrt{m_{k}M_{k}}}{M_{k}+m_{k}}+i\frac{%
\sqrt{n_{k}N_{k}}}{N_{k}+n_{k}}\right) e_{k}.  \label{3.20}
\end{equation}
\end{corollary}

The proof employs Theorem \ref{t3.2} and is similar to the one in Corollary %
\ref{c2.3}. We omit the details.

\section{Applications for Complex Numbers}

The following reverse of the generalised triangle inequality with a clear
geometric meaning may be stated.

\begin{proposition}
\label{p4.1}Let $z_{1},\dots ,z_{n}$ be complex numbers with the property
that%
\begin{equation}
0\leq \varphi _{1}\leq \arg \left( z_{k}\right) \leq \varphi _{2}<\frac{\pi 
}{2}  \label{4.1}
\end{equation}%
for each $k\in \left\{ 1,\dots ,n\right\} .$ Then we have the inequality%
\begin{equation}
\sqrt{\sin ^{2}\varphi _{1}+\cos ^{2}\varphi _{2}}\sum_{k=1}^{n}\left\vert
z_{k}\right\vert \leq \left\vert \sum_{k=1}^{n}z_{k}\right\vert .
\label{4.3}
\end{equation}%
The equality holds in (\ref{4.3}) if and only if%
\begin{equation}
\sum_{k=1}^{n}z_{k}=\left( \cos \varphi _{2}+i\sin \varphi _{1}\right)
\sum_{k=1}^{n}\left\vert z_{k}\right\vert .  \label{4.4}
\end{equation}
\end{proposition}

\begin{proof}
Let $z_{k}=a_{k}+ib_{k}.$ We may assume that $b_{k}\geq 0,$ $a_{k}>0,$ $k\in
\left\{ 1,\dots ,n\right\} ,$ since, by (\ref{4.1}), $\frac{b_{k}}{a_{k}}%
=\tan \left[ \arg \left( z_{k}\right) \right] \in \left[ 0,\frac{\pi }{2}%
\right) ,$ $k\in \left\{ 1,\dots ,n\right\} .$ By (\ref{4.1}), we obviously
have%
\begin{equation*}
0\leq \tan ^{2}\varphi _{1}\leq \frac{b_{k}^{2}}{a_{k}^{2}}\leq \tan
^{2}\varphi _{2},\qquad k\in \left\{ 1,\dots ,n\right\} 
\end{equation*}%
from where we get%
\begin{equation*}
\frac{b_{k}^{2}+a_{k}^{2}}{a_{k}^{2}}\leq \frac{1}{\cos ^{2}\varphi _{2}}%
,\qquad k\in \left\{ 1,\dots ,n\right\} ,\ \varphi _{2}\in \left( 0,\frac{%
\pi }{2}\right) 
\end{equation*}%
and%
\begin{equation*}
\frac{a_{k}^{2}+b_{k}^{2}}{a_{k}^{2}}\leq \frac{1+\tan ^{2}\varphi _{1}}{%
\tan ^{2}\varphi _{1}}=\frac{1}{\sin ^{2}\varphi _{1}},\qquad k\in \left\{
1,\dots ,n\right\} ,\ \varphi _{1}\in \left( 0,\frac{\pi }{2}\right) 
\end{equation*}%
giving the inequalities%
\begin{equation*}
\left\vert z_{k}\right\vert \cos \varphi _{2}\leq \func{Re}\left(
z_{k}\right) ,\ \ \left\vert z_{k}\right\vert \sin \varphi _{1}\leq \func{Im}%
\left( z_{k}\right) 
\end{equation*}%
for each $k\in \left\{ 1,\dots ,n\right\} .$

Now, applying Theorem \ref{t2.1} for the complex inner product $\mathbb{C}$
endowed with the inner product $\left\langle z,w\right\rangle =z\cdot \bar{w}
$ for $x_{k}=z_{k},$ $r_{1}=\cos \varphi _{2},$ $r_{2}=\sin \varphi _{1}$
and $e=1,$ we deduce the desired inequality (\ref{4.3}). The case of
equality is also obvious by Theorem \ref{t2.1} and the proposition is proven.
\end{proof}

Another result that has an obvious geometrical interpretation is the
following one.

\begin{proposition}
\label{p4.2}Let $e\in \mathbb{C}$ with $\left\vert z\right\vert =1$ and $%
\rho _{1},\rho _{2}\in \left( 0,1\right) .$ If $z_{k}\in \mathbb{C}$, $k\in
\left\{ 1,\dots ,n\right\} $ are such that%
\begin{equation}
\left\vert z_{k}-c\right\vert \leq \rho _{1},\ \ \left\vert
z_{k}-ic\right\vert \leq \rho _{2}\text{ \ for each \ }k\in \left\{ 1,\dots
,n\right\} ,  \label{4.5}
\end{equation}%
then we have the inequality%
\begin{equation}
\sqrt{2-\rho _{1}^{2}-\rho _{2}^{2}}\sum_{k=1}^{n}\left\vert
z_{k}\right\vert \leq \left\vert \sum_{k=1}^{n}z_{k}\right\vert ,
\label{4.6}
\end{equation}%
with equality if and only if%
\begin{equation}
\sum_{k=1}^{n}z_{k}=\left( \sqrt{1-\rho _{1}^{2}}+i\sqrt{1-\rho _{2}^{2}}%
\right) \left( \sum_{k=1}^{n}\left\vert z_{k}\right\vert \right) e.
\label{4.7}
\end{equation}
\end{proposition}

The proof is obvious by Corollary \ref{c2.2} applied for $H=\mathbb{C}$.

\begin{remark}
If we choose $e=1,$ and for $\rho _{1},\rho _{2}\in \left( 0,1\right) $ we
define $\bar{D}\left( 1,\rho _{1}\right) :=\left\{ z\in \mathbb{C}%
|\left\vert z-1\right\vert \leq \rho _{1}\right\} ,$ $\bar{D}\left( i,\rho
_{2}\right) :=\left\{ z\in \mathbb{C}|\left\vert z-i\right\vert \leq \rho
_{2}\right\} ,$ then obviously the intersection%
\begin{equation*}
S_{\rho _{1},\rho _{2}}:=\bar{D}\left( 1,\rho _{1}\right) \cap \bar{D}\left(
i,\rho _{2}\right) 
\end{equation*}%
is nonempty if and only if $\rho _{1}+\rho _{2}>\sqrt{2}.$

If $z_{k}\in S_{\rho _{1},\rho _{2}}$ for $k\in \left\{ 1,\dots ,n\right\} ,$
then (\ref{4.6}) holds true. The equality holds in (\ref{4.6}) if and only if%
\begin{equation*}
\sum_{k=1}^{n}z_{k}=\left( \sqrt{1-\rho _{1}^{2}}+i\sqrt{1-\rho _{2}^{2}}%
\right) \sum_{k=1}^{n}\left\vert z_{k}\right\vert .
\end{equation*}
\end{remark}

\end{document}